\documentclass[12pt, onecolumn, smallextended, runningheads,envcountsect]{article}
\usepackage{mathptmx}
\usepackage{tikz}
 \usepackage{amsfonts}
 \usepackage{amsmath}
 \usepackage{amsthm}
 \usepackage{amssymb}
 \newtheorem{thm}{Theorem}[section]
  \newtheorem{df}[thm]{Definition}
 \newtheorem{prop}[thm]{Proposition}
  \newtheorem{lmm}[thm]{Lemma}
 
 \newtheorem{rem}[thm]{Remark}
\newtheorem{ex}[thm]{Example}

\begin{document}
\begin{center}
\large{\bf{How to construct generalized van der Corput sequences  }}
\bigskip

{\large{Ingrid Carbone\footnote{University of Calabria, Department of Mathematics, Ponte P. Bucci Cubo 30B, 87036 Arcavacata di Rende (Cosenza) Italy. E-mail: i.carbone@unical.it}}}

\end{center}

\begin{abstract} 
The $LS$-sequences of points recently introduced by the author are a generalization of van der Corput sequences. They were constructed by reordering the points of the corresponding $LS$-sequences of partitions. 
Here we present another algorithm which coincides with the classical one for van der Corput sequences and is simpler to compute than the original construction. 
This algorithm is based on the representation of natural numbers in base $L+S$ and gives the van der Corput sequence in base $b$ if $L=b$ and $S=0$. 
In this construction, as well as in the van der Corput one, it is essential  the inversion of digits of the representation in base $L+S$: in this paper we also give a nice geometrical explanation of this   ``magical" operation.

\bigskip
{\bf Keywords} {Uniform distribution, sequences of partitions, van der Corput sequences, discrepancy}.

{\bf Mathematics Subject Classification (2010)} 11K06, 11K31, 11K38

\end{abstract}
\section { Introduction} 
\label{intro}

This paper contributes to enhance  the connection between the classical theory of uniformly distributed sequences of points in $[0,1[$ and the related area of uniformly distributed sequences of partitions of $[0,1[$, by showing how to reorder the points of the {\it $LS$-sequences of partitions} in order to get the {\it $LS$-sequences of points}, both introduced by the author in \cite{C}. 
If $L=b$ and $S=0$, the former generalize the sequences of partitions in base $b$, the latter generalize {\it van der Corput sequences in base $b$}. 

The construction presented here is based on the representation of positive integer numbers in base $L+S$ and coincides with the construction of van der Corput sequence in base $b$ in which, as it is well known, the representation of positive integer numbers in base $b$ is essential, as well as the inversion of these digits and the corresponding {\it radical-inverse function $\phi_b$}. 
The most interesting contribution of this paper consists in the geometrical explanation of the above mentioned inversion of digits, which appears not only for the {\it van der Corput sequences}, but (together with some more specific features) also for all the {\it $LS$-sequences of points}.

\smallskip
Let us begin by giving some definitions (for a complete overview on uniform distribution and discrepancy see \cite{KN}).

\begin{df} \label {1.1} {\rm Given a sequence  $\{Q_n\}$ of finite subsets of  $[0,1[$, with $Q_n=\{y^{(n)}_1, \dots,$ $ y^{(n)}_{t_n}\}$, 
we say that it is {\it uniformly distributed (u.d.)} if 
\begin{equation}\label{1}\lim_{n\rightarrow \infty}\frac{1}{t_n} \sum_{i=1}^{t_n} f(y_i^{(n)})= \int_0^1 f(x)\,dx\, \end{equation}
for every continuous (or Riemann integrable) function $f$ defined on $[0,1[$.}
\end{df}

\begin{df} \label {1.2} {\rm Given a sequence $\{x_n\}$ of points in  $[0,1[$, we say that it is {\it uniformly distributed (u.d.)} if the sequence of sets $Q_n=\{x_1, \dots, x_n\}$ is uniformly distributed.}
\end{df}

\begin{df} \label {1.3} {\rm Given a sequence $\{\pi_n\}$ of finite partitions of  $[0,1[$, with $\pi_n=\{[y_{i-1}^{(n)}, y_i^{(n)}[, 1\le i\le t_n\}$, $y_0^{(n)}=0$ and $y_{t_n}^{(n)}=1$, we say that it is {\it uniformly distributed (u.d.)} if the sequence of sets $Q_n=\{y^{(n)}_1, \dots, y^{(n)}_{t_n}\}$  is uniformly distributed.}
\end{df}

These definitions can be easily extended to higher dimensions, but the results of this paper concern only the one-dimensional case.

Sequences of partitions and sets have been considered rather early in the theory of uniform distribution. Partitions in base $2$ have been considered in 1935 by van der Corput  in \cite{vdC}. Knapowski \cite{Kn} studied in 1957 the sequence of partitions of $[0,1[$ into $n$ equal parts and its generalizations. Hammersley introduced  van der Corput sequences in any base $b$ and introduced a particular sequence of subsets of $I\!\!R^n$ which proved to be efficient when evaluating multidimensional integrals with Monte Carlo type methods (see \cite{DT} for these applications).
The $LS$-sequences in the multidimensional case have been studied in \cite{CIV1} and \cite{AHZ}.

More precisely, the {\it sequence of partitions in base $b$}, with $b\in \mathbb{N} $, $b\ge 2$, is defined by 
\begin{equation}\label{2}\Big\{\Big[\frac{i-1}{b^n}, \frac{i}{b^n}\Big[\ , \ 1\le i \le b^n \Big\}  \end{equation} 

\smallskip
\noindent while the corresponding sequence of points, called the {\it van der Corput sequence in base $b$}, is obtained reordering in an appropriate way the points determining them. We will not discuss this aspect in detail now, as this procedure turns out to be a special case of the algorithm we will propose in last section of this paper.
We only mention Faure's {\it scrambled van der Corput sequences} (\cite{F1, F2}) and the previous extensions to higher dimension by Hammersley \cite{Ham} and Halton \cite{Hal}.

Knapowski considered sequences of points 
\begin{equation} \label{3}\left\{\frac{1}{a_1}, \frac{2}{a_1}, \dots, \frac{a_1-1}{a_1}, \frac{1}{a_2}, \frac{2}{a_2}\dots, \frac{a_2-1}{a_2}, \frac{1}{a_n},\dots, \frac{a_n-1}{a_n}, \dots\right\}\end{equation}
related to equipartitions of $[0,1[$. He analyzed the conditions under which the sequence (\ref{3}) is u.d. (see also \cite{PS}).

Hammersley considered sets of $N$ points in $I\!\!R^d$ of the kind $(x^{(n)}_1, x^{(n)}_2, \dots x^{(n)}_d)$ for $n\le N$, where $x^{(n)}_1=\frac{n}{N}$, while for $2\le k \le d$, $x^{(n)}_k$ is the $n$-th  term of the van der Corput sequence of points in base $b_{k-1}$, the $(k-1)$-th prime number.
\smallskip

From the definitions of uniform distribution of sequences of sets, points and partitions it is clear that, when we want to evaluate the integral of $f$, the sequences of points provide a more flexible tool, as we can choose in advance any number $N$ of points $x_n$ and afterwards, if we want to improve the approximation, we can add any number of additional points using all the previously calculated values of $f$.

On the other hand, if we want to increase the number of points and we are dealing with a sequence of sets (or partitions), there are two possibilities. The worst case is when $Q_n$ and $Q_{n+k}$ have few points in common. In this case all the previously calculated values of $f$ have to be thrown away. The situation is better if $Q_n \subset Q_{n+1}$. But even this case has a drawback, as we are limited in the choice of the number of points in which to evaluate the function $f$, as 
the averages (\ref{1}) have no meaning for values between  $t_n$ and $t_{n+1}$.

\smallskip

From this observation stems the main motivation of this paper: given a sequence of u.d. sequences of partitions, it is highly desirable to associate to it a u.d. sequence of points having the lowest discrepancy.
In \cite{C} the author introduced a class of sequences of partitions (the so-called $LS$-sequences), which includes the sequence of partitions in base $b$ defined in  (\ref{2}). 
In the same article she proposed an algorithm which gives a uniformly distributed sequence of points reordering the points determining the corresponding sequence of partitions. 
When the sequence of partitions has low discrepancy, the discrepancy of the corresponding sequence of points is the same  up to a constant times a logarithmic therm (as we will recall in the next section), and this is optimal too. 
If $L \ge 2$ and $S=0$, the sequences reduce to the van der Corput sequences in base $L$.

The purpose of this paper is to present another algorithm which constructs the same sequences but in a more efficient way.

\section {\bf Preliminaries} 

Kakutani introduced in the seventies the sequence of successive $\alpha$-refinements of a partition of the unit interval \cite{K}.

\begin{df} \label {2.1} {\rm Given a finite partition $\pi$ of $[0,1[$ and given $\alpha\in \ ]0,1[$, its $\alpha$-refinement, denoted by $\alpha \pi$, is the partition obtained by subdividing all the intervals of $\pi$ having maximal length in proportion $\alpha$ and $1-\alpha$. We denote by $\alpha^n \pi$ the $\alpha$-refinement of $\alpha^{n-1} \pi$, with $\alpha^1 \pi=\alpha \pi$.}
\end{df}

He proved the following result.

\begin{thm} \label {2.2}  
 {If $\pi=\omega=\{[0,1[\}$ is the trivial partition of $[0,1[$,
the sequence $\{\alpha^n \pi \}$ is uniformly distributed.}
 \end{thm}

Kakutani's result got a considerable attention in the late seventies and early eighties, when
other authors provided different proofs of Kakutani's theorem  \cite{AF} and also
proved several stochastic versions, in which the intervals of maximal length are split according to certain probability distributions (\cite{L1}, \cite{L2}, \cite{vZ}, and \cite{BD}). 

The paper \cite{ChV} extended the notion of uniform distribution of a sequence of partitions to probability measures on complete separable metric spaces.

 In \cite{CV1} Kakutani's splitting procedure has been extended to the $n$-dimen\-sional cube with a construction which is intrinsically higher-dimensional, and in \cite{DI} it has been extended to fractals.

In \cite{CV2} the authors present a von Neumann-type theorem for sequences of partitions, showing that a sequence of partitions having an infinitesimal diameter can be reordered to a sequence of partitions which is uniformly distributed.

The paper \cite{V} introduced the following generalization of the splitting procedure.

 \begin{df} \label {2.3} \rm{For any non trivial finite partition $\rho$ of $[0,1[$, the {\it $\rho$-refinement} of a partition
$\pi$ of $[0,1[$ (denoted by $\rho \pi$) is obtained by
subdividing all the intervals of
$\pi$ having maximal length positively (or directly) homothetically to $\rho$. If   for any $n\in \mathbb{N}$ we denote by  $\rho^n \pi$ the $\rho$-refinement of $\rho^{n-1} \pi$, we get a sequence of partitions $\{\rho^n\pi\}$, called the {\it  sequence of successive $\rho$-refinements} of $\pi$.
}
\end{df}

Obviously, if $\rho =\{[0, \alpha],[\alpha, 1[\}$, then the $\rho$-refinement is just
Kakutani's $\alpha$-refinement. 

As in Kakutani's case, we can iterate the splitting procedure. The $\rho$-refinement of $\rho \pi$ will be denoted by $\rho^2 \pi$, and the meaning of $\rho^n \pi$, for $n\in I\!\!N$, is clear. The following result generalizes Theorem \ref{2.2}.

 \begin{thm} \label {2.4} {The sequence $\{\rho^n \omega\}$ is uniformly distributed.}
 \end{thm}

The recent paper by Aistleitner and Hofer \cite{AH} provides a complete solution to the question, posed in \cite{V}, for which partitions $\pi$ the sequence of partitions $\{\rho^n \pi\}$ in u.d..
\smallskip

The following result in \cite{V} connects u.d. sequences of partitions to u.d. sequences of points.

\begin{thm} \label {2.5} {\it Given a u.d. sequence of partitions $\{\pi_n\}$, a random sequential reordering of the points determining the partitions is a u.d. sequence of points with probability one.}
\end{thm}

A sequential ordering of the sequence of sets $Q_n=\{y^{(n)}_1, \dots, y^{(n)}_{t_n}\}$ is the sequence whose first $t_1$ points belong to $Q_1$, the next $t_2$ points belong to $Q_2$, and so on.     

The limit of this theorem is that it does not provide an explicit algorithm to construct a u.d. sequence of points associated to $\{\pi_n\}$.  
As we will see in the next section, for the family of $LS$-sequences introduced in \cite{C}, which are a subfamily of $rho$-refinements, this algorithm is explicitly given. 

We finish this section by giving some  very important definitions.

\begin{df} \label {2.6} {\rm Given a sequence  $\{Q_n\}$ of finite subsets of  $[0,1[$, with $Q_n=\{y^{(n)}_1, \dots,$ $ y^{(n)}_{t_n}\}$, its {\it discrepancy} is defined by the sequence
\begin{equation*} D(Q_n) = \sup_{0 \le a < b \le  1}  \left | \frac {1} {{t_n}} \sum _{j=1}^{{t_n}} \chi _{[a, \,b[} (y_j^{(n)}) - (b-a)\right |. \end{equation*}}
\end{df}

If we look at Definition 1.1, we see that discrepancy gives us the speed with which the averages $\frac{1}{t_n} \sum_{i=1}^{t_n} \chi_{[a,b[}(y_i^{(n)})$ converge to the integral of $f=\chi_{[a,b[}$. Here, as usual, $ \chi _ {[a,b[}( \cdot)$ denotes the characteristic function of the interval $[a,b[$.

\begin{df} \label {2.7}{\rm Given a sequence  $\{Q_n\}$ of finite subsets of  $[0,1[$, with $Q_n=\{y^{(n)}_1, \dots, $ $y^{(n)}_{t_n}\}$, we say that it is has {\it low discrepancy} if there exists a constant $C$ such that $t_n \  D(Q_n)\le C$.} \end{df}

 It is well known that the upper bound cannot be improved and that it is optimal (see \cite{KN}, for instance): the Knapowski sequence of sets has low discrepancy.

However, if the sequence of sets $Q_n$ is generated by a sequence of points as in Definition 1.2, it is well-known that the corresponding discrepancy cannot be so fast (see \cite{S}).

\begin{df} \label{2.8} {\rm Given a sequence  of points $\{x_n\}$ of  $[0,1[$ and $Q_n=\{x_1, \dots, x_n\}$, we say that it is has {\it low discrepancy} if there exists a constant $C$ such that $n D(Q_n)\le C \log n$.}
\end{df}

This bound is optimal and it is achieved, for example, by the van der Corput sequence in base $b$.

Let us mention that a further contribution to the theory of $\rho$-refinement has been given by Drmota and Infusino \cite{DI} who provided estimates of the discrepancy of sequences of $\rho$-refinements and, in particular, of $LS$-sequences of partitions.

\section {$LS$-sequences} 
\label{sec:3}

In this section we recall the definition of {\it $LS$-sequence of partitions} $\{ \rho_{L,S}^n \} $ and {\it $LS$-sequence of points} $\{ \xi_{L,S}^n \} $ introduced in \cite{C}.

\begin{df} \label {3.1} \rm{ Let us fix two positive integers $L\ge 1$ and $S\ge 0$ with $L+S \ge 2$ and let $\gamma \ $ be the positive solution of  the  quadratic equation $ L x + S x^2=1$ (if $S=0$ the equation is linear). Denote by $\rho_{L,S}$  the partition defined by $L$ ``long" intervals having length $\gamma$ followed by $S$ ``short" intervals having length $\gamma^2$ (if  $S=0$ all the $L$ intervals have the same length $\gamma=1/L$).
The sequence of the successive $ \rho_{L,S}$-refinements of the trivial partition $\omega$ is called $LS$-{\it sequence of partitions} and is denoted by $\{\rho^n_{L,S}\,\omega\}$ (or $\{\rho^n_{L,S}\}$ for short).
}
\end{df}

Note that for $L=b$ and $S=0$ we get the sequence of partitions  in base $b$ defined by (\ref{2}) .

If we denote with $t_n$ the total number of intervals of $ \rho^n_{L,S}$, with $l_n$ the number of  long intervals and with $s_n$ the number of  short intervals, it is very simple to see that $ t_n= l_n+s_n$, $l_n=L\,l_{n-1}+s_{n-1}$ and $s_n=S\,l_{n-1}$. We deduce that $t_n=L\, t_{n-1}+S\,t_{n-2}$ 
(with $t_0=1$ and $t_1=L+S$) which has the explicit solution

\begin{equation*} t_n= \frac{1+S\gamma} {1+S \gamma^2} \left(\frac {1} {\gamma}\right)^n- \frac{S \gamma-S \gamma^2} {1+S \gamma^2} (-S \gamma)^n .\end{equation*} 
We note that $l_n$ and $s_n$ satisfy the same difference equation, with initial conditions $l_0=1$, $l_1=L$ and $s_0=0$, $s_1=S$, respectively.

The following estimates from above and from below of the discrepancy of $LS$-sequences of partitions have been given in \cite{C}.

\begin {thm} \label{3.2}

i)  If  $S \le L$   there exist  $c_1, c'_1 >0$ such that  for any $n \in \mathbb{N}$

\begin{equation*} c'_1  \le  t_n \,D( \rho_{L,S}^n  ) \le \, c_1  \, .\end{equation*}

\bigskip
\noindent ii)  If  $S =L+1$ there exist $c_2, c'_2>0$ such that for any $n \in \mathbb{N}$

\begin{equation*}c'_2 \,{\log t_n}  \le t_n \,D( \rho_{L,S}^n  ) \le c_2 \,{\log t_n} \, .\end{equation*}

\bigskip
\noindent iii) If  $S\ge L+2$ there exist $c_3, c'_3>0$ such that for any $n \in \mathbb{N}$

  \begin{equation*} \,c'_3\, {t_n}^{1-\tau}\,  \le t_n \, D( \rho_{L,S}^n  )  \le \, c_3 \, t_n^{1-\tau}\, ,\end{equation*}

\noindent where $ 1-\tau = - \frac {\log(S \gamma)} {\log \gamma}>0$.
\end{thm}

We emphasize that each $LS$-sequence of partitions with $S\le L$ has low discrepancy. 

Let us now recall the definition given in \cite{C} of the sequence of points $ \{\xi_{L,S}^n\}$ associated to the sequence of partitions $ \{\rho_{L,S}^n\}$. 

\begin {df} \label{3.3} \rm{We fix a sequence of partitions $ \{\rho^n_{L,S} \}$. 
For each $n \ge 1$ we define the two families of functions
\begin{equation} \label {4}
{ \varphi_i^{(n+1)} (x)=x+i\gamma^{n+1}\,\,\,\,\,\,{\rm and}  \,\,\, \,\,\,\varphi_{L,j}^{(n+1)}(x)= x+L \gamma^{n+1} + j \gamma^{n+2 }  }\end{equation}

\noindent for $1 \le i \le L$ and $1 \le j \le S-1$ (if $S=0$ we have only the functions $\varphi_i^{(n+1)}$ ).  
We denote by $\Lambda_{L,S}^1$ the set made by the $t_1$ left endpoints of the intervals of $\rho_{L,S}^1\,$ ordered by magnitude. 
  If $\Lambda_{L,S}^n=\Big( \xi_1^{(n)}, \dots,  \xi_{t_n}^{(n)} \Big)$ denotes the set of the $t_n$ points defining $ \rho^n_{L,S} $, ordered in the appropriate way, the set  $\Lambda_{L,S}^{n+1}$, which consists of a reordering of the points defining $\rho^{n+1}_{L,S}$, is obtained as follows:

\noindent \begin{eqnarray}  \label {5} & &\Lambda_{L,S}^{n+1}=\Big ( \xi_1^{(n)}, \xi_2^{(n)}, \dots, \xi_{{t_n}}^{(n)},  \nonumber  \\ 
& & \varphi_1^{(n+1)}(\xi_1^{(n)} ), \dots, \varphi_1^{(n+1)} (\xi_{{l_n}}^{(n)}),  \dots,  \varphi_L^{(n+1)}(\xi_1^{(n)} ), \dots, \varphi_L^{(n+1)} (\xi_{{l_n}}^{(n)}) , \nonumber\\
& &\varphi_{L,1}^{(n+1)}(\xi_1^{(n)} ), \dots, \varphi_{L,1}^{(n+1)} (\xi_{{l_n}}^{(n)}), \dots, \varphi_{L,S-1}^{(n+1)}(\xi_1^{(n)} ), \dots, \varphi_{L,S-1}^{(n+1)} (\xi_{{l_n}}^{(n)}) \Big). \end{eqnarray}
The sequence obtained by reordering successively the points of $\rho_{L,S}^{n+1}\setminus \rho_{L,S}^{n}\,$ is called $LS$-{\it sequence of points} and is denoted by $\{\xi_{L,S}^n\} $. }
\end{df}

Let us observe that if $L=b$ and $S=0$, the sequence reduces to the van der Corput sequence in base $b$.

In \cite{C} we provided the following bounds for the discrepancy of these sequences of points.

\begin {thm} \label{3.4}

i)  If  $S \le L$   there exists  $k_1 >0$ such that for any $N \in \mathbb{N}$ 
\begin{equation*}  N \, D \Big(\xi_{L,S}^1, \xi_{L,S}^2, \dots, \xi_{L,S}^N  \Big) \le \,  k_1  {\log N}  \, .\end{equation*}

\noindent ii)  If  $S =L+1$ there exists $k_2, c'_2>0$ such that for any $N \in \mathbb{N}$ 
\begin{equation*} c'_2 \log N \le N \, D \Big( \xi_{L,S}^1, \xi_{L,S}^2, \dots, \xi_{L,S}^N  \Big) \le k_2 \,{\log^2 N} \, .\end{equation*}

\noindent iii) If  $S\ge L+2$ there exists $k_3, c'_3>0$ such that for any $N \in \mathbb{N}$
\begin{equation*} c'_3\, {N^{1-\tau}} \le N \, D \Big( \xi_{L,S}^1, \xi_{L,S}^2, \dots, \xi_{L,S}^N  \Big)  \le \, k_3 \, {N^{1-\tau}}\, {\log N} \, ,\end{equation*}
\noindent where $ 1-\tau = - \frac {\log(S \gamma)} {\log \gamma}>0$.

\end{thm}

Theorem \ref{3.2} and Theorem \ref{3.4} show that if the $LS$-sequence of partitions has low discrepancy, the corresponding $LS$-sequence of points has low discrepancy, too. Moreover,  $c'_2$ and $c'_3$ denote the same constants in both theorem as, of course,  whenever $N=t_n$ for some $n$, it follows that $D \Big( \xi_{L,S}^1, \xi_{L,S}^2, \dots, \xi_{L,S}^{N}  \Big) =D(\rho^n_{L,S})$.

Note that the lower bounds follow from Theorem \ref{3.2}. In fact  we strongly believe that  the discrepancy of all the $LS$-sequences of points coincides with the upper bound given by Theorem \ref{3.4}.

\section{Technical results}
\label{sec:4}

In this section we present an algorithm to construct $LS$-sequences of points $\{ \xi^n_{L,S}\}$, which is formally independent from the concept of sequence of partitions:  we do not need to reorder the points of each partition, as it is shown in Theorem \ref{5.3}, and therefore it is simpler to use then the algorithm introduced in \cite{C}.  Of course, the connection between each $\{ \xi^n_{L,S}\}$ and the corresponding   $\{ \rho^n_{L,S}\}$ lies behind all the procedure, and this geometric interpretation represents one of the most interesting aspects  of this new family of sequences. 

In the case $L=b$  and $S=0$, this algorithm will produce the van der Corput sequence in base $b$.
 
 From now on we fix $L,S \in \mathrm{N}$ and $\gamma$ such that $L \gamma + S \gamma^2=1$.
 
 In order to introduce  this algorithm, we represent the points of the sequence $\{ \xi_{L,S}^n \}$ writing the elements of each set $\Lambda_{L,S}^n$ defined by (\ref{5}) as images of the point $x=0$ under the compositions of suitable functions depending on $L$ and $S$ we are going to define.
 
For this purpose we introduce the functions $\psi_i$  as follows: for every $0 \le i \le L-1$ we consider
\begin{equation}\label{6} \psi_i(x)= \gamma x + i \gamma \,\,\,\,\,\,\,\,{\rm restricted\,\,\,to} \,\,\,\,\,\,\,0 \le x < 1,
\end{equation}
while for every $L \le i \le L+S-1$ and $S \ge 1$ we define
\begin{equation} \label{7}\psi_i(x)= \gamma x + L \gamma + (i-L) \gamma^2\,\,\,\,\,\,\,\,{\rm restricted\,\,\,\, to} \,\,\,\,\,\,\,0 \le x < \gamma.
\end{equation}

If $S=0$, we will consider only the functions (\ref{6}).

We observe that the functions (\ref{6}) map $[0,1[$ onto $[i \gamma\,, (i+1) \gamma\,[$ and the functions (\ref{7}) map $[0,\gamma\,[$ onto $[L \gamma+(i-L) \gamma^2, L \gamma+(i-L+1) \gamma^2[$. Therefore, the compositions   $ \psi_i\circ \psi_j $ are not defined if and only if $L \le i \le L+S-1 $ and $1 \le j \le  L+S-1$. 

Let us denote by $E_{L,S}$ the set consisting of all the pairs of indices which correspond to the "forbidden" compositions, i.e.
\begin{equation} \label{8}{E_{L,S}}=  \{ L, L+1, \dots, L+S-1\} \times \{1,\dots,L+S-1 \}.\end{equation} 

If $S=0$, the first factor is empty, so $E_{L,S}= \emptyset$.

In order to present our main result, we have to give the explicit expression of the compositions $\psi_{i_1 i_{2} \dots i_n}=\psi_{i_1} \circ\psi_{i_2}\circ \dots \psi_{i_n}$ of the functions (\ref{6}) and (\ref{7}).

\begin{lmm} \label{4.1}{For any $n \in \mathbb{N}$ and any $n$-tuple $(i_1i_2 \dots i_n)$ of elements of the set $\{0,1,\dots, L+S-1\}$ such that $(i_{h}\,,i_{h+1}) \notin E_{L,S} \,\,$ for any $1 \le h \le n-1$, we have

\begin{equation}\label{9} \psi_{i_1 i_{2} \dots i_n}(x)=\gamma^n x+ \sum_{k=1}^n {b}_k \, \gamma^k,
\end{equation} 
where $b_k=i_k\,\,$ if $ \,\,0 \le i_k\le L-1$, while
$\,b_k=L+(i_k-L)\gamma \,\,$ if $\,\,L \le i_k \le L+S-1$. 
}
\end{lmm}
\proof We prove (\ref{9}) by induction on $n$. Of course, we need to consider only the case $S \ge 1$.

If $n=2$, because of (\ref{6}) and (\ref{7}) we have only three kinds of compositions $\psi_{i_1, i_2}$: 

\smallskip
\noindent a) if $0 \le i_1, i_2 \le L-1$, we have 
$\psi_{i_1, i_2}(x)= 
\gamma^2x+ \sum_{k=1}^2 b_k \,\gamma^k , $
 where $b_k=i_k$; 

\smallskip
\noindent b) if $0 \le i_1 \le L-1$ and $L \le i_2 \le L+S-1$, a simple calculation gives (\ref{9}), with $b_1=i_1$ and $b_2=L+(i_2-L) \gamma$;

\smallskip
\noindent c) if $L \le i_i \le L+S-1$ and $ 0 \le i_2 \le L-1$, we obtain the same expression, where $b_1=L+(i_1-L)\gamma$ and $b_2=i_2 $. 

\smallskip
Suppose now that (\ref {11}) holds for $n \ge 1$ and let us prove it for $n+1$, namely that
\begin{equation}\label{10} \psi_{i_1i_2 \dots i_n i_{n+1}}(x)=\gamma^{n+1}x+ \sum_{k=1}^{n+1}b_k \, \gamma^k
\end{equation}
for all $i_{n+1}\neq 0$. Of course, at this stage we have to distinguish between the two kinds of functions (\ref{6}) and (\ref{7}).

If $0 \le i_{n+1} \le L-1$, we simply obtain (\ref{10}) with $b_{n+1}=i_{n+1}$.

If $L \le i_{n+1} \le L+S-1$, simple calculations give again (\ref{10}), with $b_{n+1}=L+(i_{n+1}-L)\gamma$. So the induction is complete. $\qed$

\begin{rem} We point out that if $S=0$, the coefficients $b_k$ in formula (\ref{9}) are all equal to $i_k$ as $E_{L,S}= \emptyset$.
\end{rem}

As a consequence of the previous lemma we obtain the following result.

\begin{prop}\label{4.3}   The first $t_n$ points of the $LS$-sequence  $\{ \xi_{L,S}^n \}$ are
\begin{equation}\label{11}\Lambda_{L,S}^n= \{\psi_{i_1 i_{2} \dots i_n} (0) \}, \end{equation} 
where $(i_{h}\,,i_{h+1}) \notin E_{L,S}\, $ for any $1 \le h \le n-1$ and $(i_n, i_{n-1},\dots,i_1)$ are the $n$-tuples of elements of the set $\{0,1,\dots, L+S-1\}$, ordered with respect to the magnitude of the numbers $i_n i_{n-1} \dots i_1$ in base $L+S$. 
\end{prop}

\proof 

Taking (\ref{6}) and (\ref{7}) into account, we see that 
\begin{equation*}\Lambda_{L,S}^1=\{\psi_{0}(0), \psi_{1}(0), \dots,\psi_{L+S-1}(0) \},
\end{equation*}
which are exactly the $L+S$ left endpoints of the intervals of $\rho_{L,S}^1$ ordered by magnitude, so (\ref{11}) is true for $n=1$.

Fix any integer $n \ge 1$ and suppose (\ref{11}) is true for $n$. 

If $i_{n+1}=0$, as $\psi_0(0)=0$ we have 
\begin{equation*}\psi_{i_1i_2 \dots i_n i_{n+1}}(0)=\psi_{i_1i_2 \dots i_n}(0),
\end{equation*}
and we simply observe that this way we get all the points of $\Lambda_{L,S}^n$. All the new points of $\Lambda_{L,S}^{n+1} \setminus \Lambda_{L,S}^n$ are obtained when $i_{n+1} \not=0$ as follows.

If $1 \le i_{n+1} \le L$ we have
\begin{eqnarray*} \psi_{i_1i_2 \dots i_n i_{n+1}}(0)&=&
i_{n+1}\gamma^{n+1}+\sum_{k=1}^{n}b_k \, \gamma^k=\psi_{i_1i_2 \dots i_n}(0)+ i_{n+1}\gamma^{n+1}\\
&=&\varphi_{i_{n+1}}^{(n+1)}(\psi_{i_1i_2 \dots i_n}(0)),
\end{eqnarray*}
according to (\ref{4}). 

If $L+1 \le i_{n+1} \le L+S-1$, we set $j_{n+1}=i_{n+1}-L$, with $1 \le j_{n+1} \le S-1$, and write 
 \begin{eqnarray*} &&\psi_{i_1i_2 \dots i_n i_{n+1}}(0)=
 (L+j_{n+1}\gamma) \gamma^{n+1}+\sum_{k=1}^{n}b_k \, \gamma^k\\
 &=&\psi_{i_1i_2 \dots i_n}(0)+ L \gamma^{n+1}+j_{n+1}\gamma^{n+2}=\varphi_{L, j_{n+1}}^{(n+1)}(\psi_{i_1i_2 \dots i_n}(0)).
\end{eqnarray*}

\smallskip
At this point we recall that the set $\Lambda_{L,S}^{n+1}\setminus \Lambda_{L,S}^n$ is made by $(L+S-1)l_n$ points as we evaluate the $L+S-1$ functions defined by (\ref{4}) at the first $l_n$ points of $\Lambda_{L,S}^n$ (see (\ref{5})). 

It remains only to prove that the total number of the points of $\Lambda_{L,S}^{n+1}\setminus \Lambda_{L,S}^{n}$  we add in the way described above is exactly $(L+S-1)l_n$. We denote by $d_{L,S}^{(n)}$ the total number of $n$-tuples containing two consecutive digits of the set $E_{L,S}$ defined by (\ref{8}). Of course, we have $d_{L,S}^{(n)}=(L+S)^n-t_n$.

We observe, in fact, that the points of $\Lambda_{L,S}^{n+1}\setminus \Lambda_{L,S}^{n}$ come from $L+S-1$ blocks of all the  $(n+1)$-tuples of elements of the set $\{0,1,\dots, L+S-1\}$, which are $(L+S)^{n}$, ordered by magnitude of the numbers in base $L+S$.  
Moreover, the total number of not admissible $(n+1)$-tuples which are contained in these blocks is given by the difference between $(L+S-1)(L+S)^{n}$ and the total number of forbidden ones, which is exactly  $d_{L,S}^{(n+1)}-d_{L,S}^{(n)}$. 
In other words, 
\begin{eqnarray*} &&(L+S-1) (L+S)^{n}-\Big(d_{L,S}^{(n+1)}-d_{L,S}^{(n)}\Big)\\
&&=(L+S-1) (L+S)^{n}-\big((L+S)^{n+1}-t_{n+1} -(L+S-1)^{n}+t_{n}\big)  \\
&&= t_{n+1}-t_{n}=(L+S-1)l_{n}.
\end{eqnarray*}

The proposition is completely proved. $\qed$ 

\section{The new construction}
\label{sec:5}

We are now ready to give the construction of $LS$-sequences of points, which is deeply connected to the explicit geometric representation given in Proposition \ref{4.3}. More precisely, we want to show how to get the points of these sequences starting from  the representation  of natural numbers in base $L+S$. This  algorithm  is very simple and allows us to compute directly the points of each $\{ \xi_{L, S}^n\}$. 

This algorithm is, in fact, a generalization of the algorithm which produces the van der Corput sequences. 

Any natural number $n\ge 1$ can be expressed in base $b\ge 1$ as
\begin {equation}\label{12} n=\sum_{k=0}^M \, a_k(n) \,b^k \,,\end{equation} 

\noindent with $a_k(n) \in \{0,1, \dots, b-1  \}$ for all $0 \le k \le M$, with $M={\cal b}{\log_b n} \cal c$ (here and in the sequel ${\cal b} \cdot \cal c$ denotes the integer part).

The expression (\ref{12}) leads to the {\it representation in base $b$}  
\begin{equation}\label{13}[n]_b=a_M(n) a_{M-1}(n)\dots a_0(n)\,.\end{equation}
If $n=0$, we write $[0]_b=0$.

\smallskip
 The representation of $n$ in base $b$ given by (\ref{12}) is used to define the  {\it radical-inverse function $\phi_b$} on $\mathbb{N}$ which associates to the string of digits (\ref{13}) the number 
 \begin{equation}\label{14}\phi_b(n)=\sum_{k=0}^M a_k(n) b^{-k-1}  \,,\end{equation} 
whose binary representation is  $0.a_0(n) a_1(n) \dots a_M(n)$.

Of course $0 \le \phi_b(n) <1$ for all $n \ge 0$ and the $b$-radix notation of $\phi_b(n)$ actually is $[\phi_b(n)]_b=0.a_0(n) a_1(n) \dots a_M(n)$. 

The sequence $\{ \phi_b(n)\}_{n\ge 0}$ is the {\it van der Corput sequence in base $b$}. 

 \begin{df}\label{5.1} {\rm  We denote by $ \mathbb{N}_{L,S} $ the set of all positive integers $n$,  ordered by magnitude, with  $[n]_{L+S}= a_M(n)\,a_{M-1}(n)\, \dots a_0(n)$ such that  $(a_{k}(n),a_{k+1}(n))\notin E_{L,S}\,\,$ for all $0 \le k\le M-1$.
 If $S=0$, we have $ \mathbb{N}_{L,S} = \mathbb{N}$. } \end{df}

\smallskip
 \begin{df}\label{5.2} {\rm For all $n\in \mathbb{N}_{L,S}$  we define the} $LS$-radical inverse function {\rm as follows:
 \begin{equation}\label{15} \phi_{L,S}(n)=\sum_{k=0}^M \,  \tilde a_k(n) \, \gamma ^{k+1} \,,\end{equation} 
where   $ \tilde a_k(n)= a_k(n) $ if $0 \le a_k(n) \le L-1 $ and  
$\tilde a_k(n) =L+\gamma(a_k(n)-L)$ if $ L\le a_k(n) \le L+S-1$. 
(If $S=0$, (\ref{15}) coincides with the radical inverse function (\ref{14}).)} \end{df}

\smallskip
Taking the above definitions into account, we are able to present the main result of this paper.

 \begin{thm}\label{5.3} Any $LS$-sequence of points coincides with the sequence $\{\phi_{L,S}(n)\}$ defined on $\mathbb{N}_{L,S}$.
 
\proof Fix $n \in \mathbb{N}_{L,S}$. Taking Lemma \ref{4.1} and Definition \ref{5.2}  into account, we write 
 \begin{equation*} \phi_{L,S}(n)=\psi_{n_0 n_1  \dots  n_M}(0), 
 \end{equation*} 
 which proves the theorem. $\qed$ 
\end{thm} 

\begin{rem} \label{5.4} If $L=b$ and $S=0$, we have $\phi_{b,0}=\phi_{b}$ and, therefore, the $b,0$-sequence is the van der Corput sequence in base $b$. 
All the results contained in this section can be re-written in this particular case, a.
\end{rem}

\begin{ex}\label{5.5} {\rm In this example we discuss the $1,1$-sequence of points in order to better illustrate the algorithm described in this section and to underline analogies and differences with the van der Corput sequence in base $2$.
For yet another apporach to the $1,1$-sequence of points, see \cite{CIV2}.

We consider the  Kakutani-Fibonacci sequence of points $\{ \xi_{1,1}^n\}$ corresponding to $ L=S=1$ and $\gamma={1\over 2}(\sqrt{5}-1)$. See Fig. 1. 

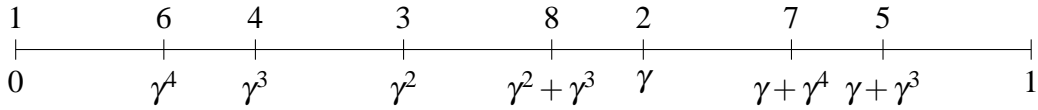
\begin{figure}[h!]

\begin{center}
\begin{tikzpicture}[scale=13.5]
\draw (0, 0) -- (1,0);
\draw (0,-0.01) node[below, black]{0} -- (0,0.01) node[above, black]{1};
\draw (1,-0.01) node[below, black]{1} -- (1,0.01);
\draw (0.61803,-0.01) node[below, black]{$\gamma$}-- (0.61803,0.01) node[above, black]{2};
\draw (0.38196,-0.01) node[below, black]{$\gamma^{2}$}-- (0.38196,0.01) node[above, black]{3};
\draw (0.23606,-0.01) node[below, black]{$\gamma^{3}$}-- (0.23606,0.01) node[above, black]{4};
\draw (0.85410,-0.01) node[below, black]{$\gamma +\gamma^{3}$}-- (0.85410,0.01) node[above, black]{5};
\draw (0.14589,-0.01) node[below, black]{$\gamma^{4}$}-- (0.14589,0.01) node[above, black]{6};
\draw (0.76393,-0.01) node[below, black]{$\gamma +\gamma^{4}$}-- (0.76393,0.01) node[above, black]{7};
\draw (0.52786,-0.01) node[below, black]{$\gamma^2 +\gamma^{3}$}-- (0.52786,0.01) node[above, black]{8};
\end{tikzpicture}
\end{center}
\caption{The first $8$ points of $\{ \xi_{1,1}^n\}$.
}
\end{figure}

\begin{figure}[h!]
\begin{center}
\begin{tikzpicture}[scale=7]
\draw (0,0) node[below]{0} -- (1,0) node[below]{1} -- (1,1) -- (0,1)node[left]{1} -- (0,0);
\draw [thick] (0,0.61803) node[left]{$\gamma$} -- (0.61803,1)
node[midway, sloped,above] {$\psi_1$};
\draw [dashed] (0.61803,0) node[below]{$\gamma$} -- (0.61803,1);
\draw [dashed] (0,0.61803) -- (1,0.61803);
\draw [thick] (0,0) -- (1,0.61803)
node[midway, sloped,above] {$\psi_0$};
\end{tikzpicture}
\end{center}
\caption{The functions $\psi_0$ and $\psi_1$ associated to $\{ \xi_{1,1}^n\}$}
\end{figure}
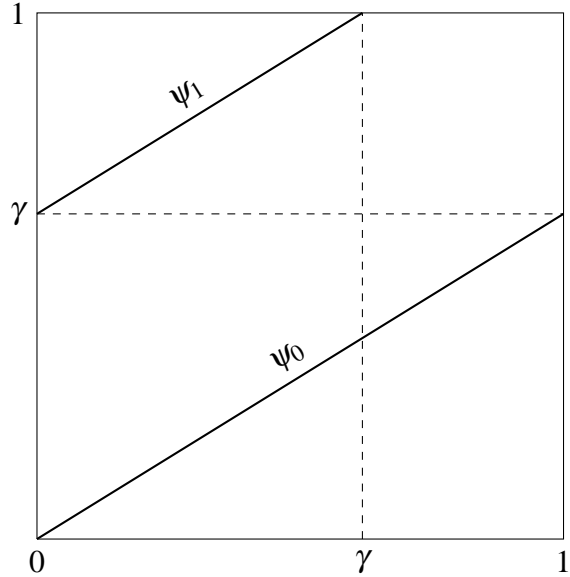

The functions defined by (\ref {8}) and (\ref{7}) reduce to  $\psi_0(x)=\gamma x$ for  $0\le x < 1$ and $\psi_1(x)= \gamma \,x + \gamma$  for $0 \le x < \gamma$. See the Fig. 2.

According to Definition \ref{5.2}, $ \mathbb{N}_{1,1}$ is the set of all natural numbers $n$ such that the binary representation  (\ref{13}) does not contain two consecutive digits equal to $1$. Moreover, the $(1,1)$-radical inverse function defined by (\ref{12}) on $ \mathbb{N}_{1,1}$ is
\begin{equation*} \phi_{1,1}(n)=\sum_{k=0}^M \, a_k(n) \, \gamma ^{k+1} \,,\end{equation*}

\noindent with the same coefficients $a_k(n)$ of the binary representation $n$ given by (\ref{13}) for $b=2$.

By Theorem \ref{5.3}, the Kakutani-Fibonacci sequence of points $\{\xi_{1,1}^n\}$ coincides with the sequence $\{\phi_{1,1}(n) \}$.

The following table shows the construction of the first $12$ points of $\{ \xi_{1,1}^n\}$, where the first column contains the first $12$ elements of $\mathbb{N}_{1,1}$.

\begin{center}

\begin{tabular}{ccccccccccc}

$0\rightarrow$ & $[0]_2$ & $=$ & $0$ $\,\,\,\,\,\,\,\,\,\,\,\,\, \rightarrow$ & $0.0$ & $=$ &
$[\phi_{1,1}(0)]_2$ & $\rightarrow$ &$0$ &$=\xi_{1,1}^1$& \\

$1\rightarrow$ & $[1]_2$ & $=$ & $1$ $\,\,\,\,\,\,\,\,\,\,\,\,\rightarrow$ & $0.1$ & $=$ &
$[\phi_{1,1}(1)]_2$ & $\rightarrow$ & $\gamma$ &$=\xi_{1,1}^2$&\\

$2\rightarrow$ & $[2]_2$ &$=$ & $10$ $\,\,\,\,\,\,\,\,\,\rightarrow$ & $0.01$ &$=$ &
$[\phi_{1,1}(2)]_2$ & $\rightarrow$ & $\gamma^2$ &$=\xi_{1,1}^3$&\\

$4\rightarrow$ & $[4]_2$ & = & $100$ $\,\,\,\,\,\,\rightarrow$ & $0.001$ & $=$ &
$[\phi_{1,1}(4)]_2$ & $\rightarrow$ & $\gamma^{3}$ &$=\xi_{1,1}^4$&\\

$5\rightarrow$ & $[5]_2$& = & $101$ $\,\,\,\,\,\,\,\rightarrow$ & $0.101$ & $=$ &
$[\phi_{1,1}(5)]_2$ & $\rightarrow$ & $\gamma+\gamma^{3}$&$=\xi_{1,1}^5$&\\

$8\rightarrow$ & $[8]_2$ & = & $1000$ $\,\,\,\,\rightarrow$ & $0.0001$ & $=$
&$[\phi_{1,1}(8)]_2$ & $\rightarrow$ & $\gamma^{4}$ &$=\xi_{1,1}^6$&\\

$9\rightarrow$ & $[9]_2$ & = & $1001$ $\,\,\,\,\rightarrow$ & $0.1001$ & $=$ &
$[\phi_{1,1}(9)]_2$ & $\rightarrow$ & $\gamma+\gamma^{4}$ &$=\xi_{1,1}^7$&\\

$10\rightarrow$ & $[10]_2$ & = & $1010$ $\,\,\,\,\rightarrow$ & $0.0101$& $=$
& $[\phi_{1,1}(10)]_2$ & $\rightarrow$ & $\gamma^2+\gamma^{4}$&$=\xi_{1,1}^8$&\\
$16\rightarrow$ & $[16]_2$ &= & $10000$ $\rightarrow$ & $0.00001$& $=$ &
$[\phi_{1,1}(16)]_2$ & $\rightarrow$ & $\gamma^{5}$&$=\xi_{1,1}^9$&\\

$17\rightarrow$ & $[17]_2$ & = & $10001$ $\rightarrow$ & $0.10001$ & $=$ &
$[\phi_{1,1}(17)]_2$ & $\rightarrow$ & $\gamma+\gamma^{5}$&$=\xi_{1,1}^{10}$& \\

$18\rightarrow$ & $[18]_2$ & = & $10010$ $\rightarrow$ & $0.01001$ & $=$
&$[\phi_{1,1}(18)]_2$ & $\rightarrow$ & $\gamma^2+\gamma^{5}$ &$=\xi_{1,1}^{11}$&\\

$20\rightarrow$ & $[20]_2$ & = & $10100$ $\rightarrow$ & $0.00101$& $=$
& $[\phi_{1,1}(20)]_2$ & $\rightarrow$ & $\gamma^3+\gamma^{5}$&$=\xi_{1,1}^{12}$&\\
\end{tabular}

\end{center}

}
\end{ex}

\bigskip

\noindent {\bf Acknowledgements}

The author wishes to  express her gratitude to Aljo\v{s}a Vol\v{c}i\v{c} for his critical reading of the manuscript and his useful suggestions.


\begin{thebibliography}{}

\bibitem {AF} R.L. Adler,  L. Flatto, Uniform distribution of Kakutani's interval splitting procedure. {\it Z. Wahrscheinlichkeitsthorie verw. Gebiete} {\bf 38} (1977), 253-259.

\bibitem {AH} C. Aistleitner, M. Hofer, Uniform distribution of generalized Kakutani's sequences of partitions. {\it Ann. Mat. Pura Appl.}, DOI: 10.1007/s10231-011-0235-9.

\bibitem{AHZ}C. Aistleitner, M. Hofer, V. Ziegler, On the uniform distribution modulo 1 of multidimensional LS-sequences. {\it Ann. Mat. Pura Appl.}, to appear.

\bibitem {BD} M.D. Brennan,  R. Durrett, Splitting intervals, {\it Ann. Probability} {\bf 14} no. 3 (1986), 1024-1036.

 \bibitem  {C} I. Carbone,  Discrepancy of $LS$ sequences of partitions and points,  {\it Ann. Mat. Pura Appl. (4)} {\bf 191} (2012) no. 4, 819Ð844.
 
\bibitem{CIV1} I. Carbone, M.R. Iac\`o, A. Vol\v{c}i\v{c}, $LS$-sequences of points in the unit square.  	 {\it submitted }(2012) arXiv:1211.2941.
 
 \bibitem {CIV2} I. Carbone, M. R. Iac\`{o}, A. Vol\v{c}i\v{c}, A dynamical system approach to the Kakutani-Fibonacci sequence. {\it Ergod. Th. and Dynam. Sys.} (2013). doi:10.1017/etds.2013.20. 
   
\bibitem  {CV1} I. Carbone, A. Vol\v{c}i\v{c}, Kakutani splitting procedure in
higher dimension, {\it  Rend. Ist. Matem. Univ. Trieste} {\bf 39}  (2007), 119-126.

\bibitem {CV2} I. Carbone, A. Vol\v{c}i\v{c}, A von Neumann theorem for uniformly distributed sequences of partitions, {\it Rend. Circ. Mat.
Palermo }{\bf 60} n. 1-2 (2011), 83-88.

 \bibitem  {ChV} F. Chersi,  A. Vol\v{c}i\v{c},  $\lambda$-equidistributed sequences of partitions and a theorem of the de Bruijn-Post type,  {\it Ann. Mat. Pura Appl. (4) {\bf 162}} (1992), 23-32.


\bibitem {DI} M. Drmota, M. Infusino,   On the discrepancy of some generalized Kakutani's sequences of partitions, {\it  Unif. Distrib. Theory} {\bf 7} (2012) no. 1, 75-104.  

\bibitem {DT} M. Drmota,  R. F. Tichy,  {\it Sequences, discrepancies and applications},  Lecture Notes in Mathematics {\bf 1651}, Springer Verlag, Berlin, 1997. 


\bibitem {F1} H. Faure, Discr\'epances de suites associ\'es \`a un syst\`eme de num\'eration (en dimension un),  {\it Bull. Soc. Math. France} {\bf 109} (1981) no. 2, 143-182.

\bibitem {F2} H. Faure,  Good permutations for Extreme Discrepancy,  {\it J. Number Theory}  {\bf 42} (1992), 47-56.

\bibitem {Hal} J. H. Halton, On the efficiency of certain quasi-random sequences of points in evaluating multi-dimensional integrals, {\it Numer. Math.} {\bf 2} (1960), 84-90.

\bibitem {Ham} J. M. Hammersley, Monte-Carlo methods for solving multivariate problems, {\it Ann. New York Acad. Sci.} {\bf 86} (1960), 844-874.

\bibitem {IV} M. Infusino, A. Vol\v{c}i\v{c}, Uniform distribution on fractals, {\it Unif. Distrib. Theory} {\bf 4} (2009), n. 2, 47-58.  

\bibitem {K} S. Kakutani, A problem on equidistribution on the unit
interval $[0,1[$, {\it  Measure theory (Proc. Conf., Oberwolfach, 1975)},  pp. 369-375. {\it Lecture Notes in Math.} {\bf 541}, Springer, Berlin, 1976.

\bibitem {Kn} S. Knapowski, Uber ein Problem der Gleichverteilungstheorie. {\it Colloq. Math.} {\bf 5} (1958), 8-10.

\bibitem {KN} L. Kuipers, H. Niderreiter, {\it Uniform distribution of sequences.  Pure and Applied Matematics}. Wiley-Interscience, New York-London-Sidney, 1974.

\bibitem {L1} J. C. Lootgieter,  Sur la r\'epartition des suites de Kakutani. I. {\it  Ann. Inst. H. Poincar\'e Sect. B (N.S.)} {\bf  13} no. 4 (1977),  385-410.

\bibitem {L2} J. C. Lootgieter, Sur la r\'epartition des suites de Kakutani. II. {\it  Ann. Inst. H. Poincar\'e Sect. B (N.S.)} {\bf  14} no. 3 (1978),   279-302. 

\bibitem {PS} \v{S}. Porubsk\'y, T. \v{S}al\'at, O. Strauch, On a class of uniformly distributed sequences. {\it Math. Slovaca} {\bf 40} (1990), no. 2, 143-170.

\bibitem {S} W. M. Schmidt, Irregularities of distribution. VII, {\it Acta Arith.} {\bf 21} (1972), 45-50.

\bibitem {vdC} J. C. van der Corput, Verteilungsfunktionen, {\it Proc. Ned. Akad. v. Wet.} {\bf 38} (1935), 813-821.

\bibitem {vZ} W. R. van Zwet, A proof of Kakutani's conjecture on random subdivision of longest intervals. {\it Ann. Probability} {\bf  6} no. 1 (1978),  133-137

\bibitem {V} A. Vol\v{c}i\v{c},  A generalization of Kakutani's splitting procedure. {\it Ann. Mat. Pura Appl.}, Nuova Serie  {\bf 190}  (2011), no. 1, 45-54.



\end{thebibliography}
\end{document}